# ON MULTIDIMENTIONAL PYTHAGOREAN NUMBERS

D.A. Sardelis –T.M Valahas

## I. OVERVIEW

To represent positive integers by regular patterns on a plane or in three-dimensional space remains an active trend in number theory (e.g. see [**3**]) which may be traced back (see [**1**]) to the works of its great pioneers − Gauss, Euler, Fermat, etc. − and even further back, to the founders of this conception, the Pythagoreans. The aim of the present article is to explore the possibility of extending the representation framework for integers to spaces with more than three dimensions. Thus, taking up a definition of polygonal numbers given by Diophantus [**2**] and also by Nicomachus [**4**] and generalizing the Pythagorean concept of gnomon, one is led through quite elementary means to a single, unified definition of multidimensional number formations henceforth called **hypersolids**. Viewing these numbers from different perspectives, several intrinsic symmetries become manifest which are worthy of further exploration as they may be of use in number theory and/or other fields.

## II. PLANE POLYGONAL NUMBERS

Polygonal numbers may be defined as follows: *Given an arithmetical progression with the first term 1 and common difference d, the sum of n terms is the nth polygonal with* $k = d + 2$ *vertices and k sides with n units each.*
Let $x_r$ denote the *r*th term of the given progression and $p(d,n)$ the *n*th polygonal number. Then, $x_r = 1+(r-1)d$ and $p(d,n)$ can be written

$$p(d,n) = \sum_{r=1}^{n} x_r = \frac{1}{2}n[2+(n-1)d]. \qquad (1)$$

Since the $x_n$ terms generate the polygonal number sequence in the sense that

$$p(d,n) = p(d,n-1) + x_n, \qquad (2)$$

they were called the gnomons of polygonals and henceforth they will be referred to as *n*-gnomons.
Similarly, one may also define as *d*-gnomons the differences of any two polygonals with the same *n* and consecutive *d* : $p(d,n) - p(d-1,n)$ . Thus, using (1), we have

$$p(d,n) = p(d-1,n) + p(1,n-1). \qquad (3)$$

Therefore, the *d*-gnomons of polygonals are the triangular numbers.
Table 1 displays the polygonal numbers for $d,n \le 10$ and their corresponding gnomons −*d* and *n*− denoted as (*d*) and (*n*), respectively (plane numbers are graphically illustrated in Appendix A).

**Table 1 Polygonal Numbers and Their Gnomons [(n),(d)]**

| d \ n | 1 | 2 | 3 | 4 | 5 | 6 | 7 | 8 | 9 | 10 | (n) |
|---|---|---|---|---|---|---|---|---|---|---|---|
| **1** | 1 | 3 | 6 | 10 | 15 | 21 | 28 | 36 | 45 | 55 | n |
| **2** | 1 | 4 | 9 | 16 | 25 | 36 | 49 | 64 | 81 | 100 | 2n–1 |
| **3** | 1 | 5 | 12 | 22 | 35 | 51 | 70 | 92 | 117 | 145 | 3n–2 |
| **4** | 1 | 6 | 15 | 28 | 45 | 66 | 91 | 120 | 153 | 190 | 4n–3 |
| **5** | 1 | 7 | 18 | 34 | 55 | 81 | 112 | 148 | 189 | 235 | 5n–4 |
| **6** | 1 | 8 | 21 | 40 | 65 | 96 | 133 | 176 | 225 | 280 | 6n–5 |
| **7** | 1 | 9 | 24 | 46 | 75 | 111 | 154 | 204 | 261 | 325 | 7n–6 |
| **8** | 1 | 10 | 27 | 52 | 85 | 126 | 175 | 232 | 297 | 370 | 8n–7 |
| **9** | 1 | 11 | 30 | 58 | 95 | 141 | 196 | 260 | 333 | 415 | 9n–8 |
| **10** | 1 | 12 | 33 | 64 | 105 | 156 | 217 | 288 | 369 | 460 | 10n–9 |
| **(d)** | 0 | 1 | 3 | 6 | 10 | 15 | 21 | 28 | 36 | 45 | |

## III. SOLID PYRAMIDAL NUMBERS

Just as polygonal numbers are produced by summing arithmetical progressions, so solid pyramidal numbers are obtained by piling successive polygonals one upon the other.

Solid pyramidals may be defined as follows: *Given a polygonal sequence* $p(d,1), p(d,2), \ldots,$ *the sum of n terms is the nth pyramidal number with k+1 = d+3 vertices and k+1 edges of n units each.* Evidently, the base of the so-defined solid number is the *n*th polygonal, hence its name as triangular, square, pentagonal, etc., pyramid.

Let *P(d,n)* denote the *n*th pyramidal number. Then, applying (1) and the summation formulas

$$\sum_{r=1}^{n} r = \frac{1}{2} n(n+1), \quad \sum_{r=1}^{n} r^2 = \frac{1}{6} n(n+1)(2n+1), \tag{4}$$

*P(d,n)* can be written

$$P(d,n) = \sum_{r=1}^{n} p(d,r) = \frac{1}{6} n(n+1)[3+(n-1)d]. \tag{5}$$

It readily follows that pyramidals have as *n*-gnomons the corresponding polygonals, i.e.,

$$P(d,n) = P(d,n-1) + p(d,n), \tag{6}$$

and as *d*-gnomons, the triangular pyramidals, i.e.,

$$P(d,n) = P(d-1,n) + P(1,n-1). \tag{7}$$

Table 2 displays the pyramidal numbers for $d, n \le 10$ and their *(d)* and *(n)* gnomons (solid numbers are graphically illustrated in Appendix B).

**Table 2 Pyramidal Numbers and Their Gnomons [(n), (d)]**

| d \ n | 1 | 2 | 3 | 4 | 5 | 6 | 7 | 8 | 9 | 10 | (n) |
|---|---|---|---|---|---|---|---|---|---|---|---|
| **1** | 1 | 4 | 10 | 20 | 35 | 56 | 84 | 120 | 165 | 220 | n(n+1)/2 |
| **2** | 1 | 5 | 14 | 30 | 55 | 91 | 140 | 204 | 285 | 385 | $n^2$ |
| **3** | 1 | 6 | 18 | 40 | 75 | 126 | 196 | 288 | 405 | 550 | n(3n–1)/2 |
| **4** | 1 | 7 | 22 | 50 | 95 | 161 | 252 | 372 | 525 | 715 | n(2n–1) |
| **5** | 1 | 8 | 26 | 60 | 115 | 196 | 308 | 456 | 645 | 880 | n(5n–3)/2 |
| **6** | 1 | 9 | 30 | 70 | 135 | 231 | 364 | 540 | 765 | 1045 | n(3n–2) |
| **7** | 1 | 10 | 34 | 80 | 155 | 266 | 420 | 624 | 885 | 1210 | n(7n–5)/2 |
| **8** | 1 | 11 | 38 | 90 | 175 | 301 | 476 | 708 | 1005 | 1375 | n(4n–3) |
| **9** | 1 | 12 | 42 | 100 | 195 | 336 | 532 | 792 | 1125 | 1540 | n(9n–7)/2 |
| **10** | 1 | 13 | 46 | 110 | 215 | 371 | 588 | 876 | 1245 | 1705 | n(5n–4) |
| **(d)** | 0 | 1 | 4 | 10 | 20 | 35 | 56 | 84 | 120 | 165 | |

## IV. FOUR-DIMENSIONAL SOLID NUMBERS

Despite religious and philosophical beliefs, there are no logical grounds for confining intellectual operations to the three-dimensional space of sense perception. Subsequently, the procedure employed above can be extended to generate solid numbers in higher dimensions. Thus, by piling successive pyramidals one obtains four-dimensional solid numbers.

Four-dimensional pyramidal numbers (hypersolids) may be defined as follows: *Given a pyramidal sequence* $P(d,1), P(d,2), \ldots,$ *the sum of n terms is the nth hypersolid in four dimensions with* $k+2 = d+4$ *vertices and k + 2 edges of n units each.* The "base" of the so-defined hypersolid is the *n*th pyramidal number.
Let $\Pi(d,n)$ denote the *n*th hypersolid. Then, using (4), (5) and the summation formula

$$\sum_{r=1}^{n} r^3 = \left[\frac{n(n+1)}{2}\right]^2, \tag{8}$$

$\Pi(d,n)$ takes the form

$$\Pi(d,n) = \sum_{r=1}^{n} P(d,r) = \frac{1}{24} n(n+1)(n+2)[4+(n-1)d]. \tag{9}$$

It follows that four-dimensional solids have as *n*-gnomons the corresponding pyramidal numbers, i.e.,

$$\Pi(d,n) = \Pi(d,n-1) + P(d,n), \tag{10}$$

and as *d*-gnomons, the lowest order four-dimensional solid numbers, i.e., the pentahedrals $\Pi(1,n-1)$:

$$\Pi(d,n) = \Pi(d-1,n) + \Pi(1,n-1). \tag{11}$$

Table 3 displays the four-dimensional solid numbers for $d,n \leq 10$ and their (*d*) and (*n*) gnomons.

**Table 3 Four-Dimensional Solid Numbers and Their Gnomons [(n), (d)]**

| d \ n | 1 | 2 | 3 | 4 | 5 | 6 | 7 | 8 | 9 | 10 | (n) |
|---|---|---|---|---|---|---|---|---|---|---|---|
| **1** | 1 | 5 | 15 | 35 | 70 | 126 | 210 | 330 | 495 | 715 | n(n+1)(n+2)/6 |
| **2** | 1 | 6 | 20 | 50 | 105 | 196 | 336 | 540 | 825 | 1210 | n(n+1)(2n+1)/6 |
| **3** | 1 | 7 | 25 | 65 | 140 | 266 | 462 | 750 | 1155 | 1705 | $n^2$(n+1)/2 |
| **4** | 1 | 8 | 30 | 80 | 175 | 336 | 588 | 960 | 1485 | 2200 | n(n+1)(4n–1)/6 |
| **5** | 1 | 9 | 35 | 95 | 210 | 406 | 714 | 1170 | 1815 | 2695 | n(n+1)(5n–2)/6 |
| **6** | 1 | 10 | 40 | 110 | 245 | 476 | 840 | 1380 | 2145 | 3190 | n(n+1)(2n–1)/2 |
| **7** | 1 | 11 | 45 | 125 | 280 | 546 | 966 | 1590 | 2475 | 3685 | n(n+1)(7n–4)/6 |
| **8** | 1 | 12 | 50 | 140 | 315 | 616 | 1092 | 1800 | 2805 | 4180 | n(n+1)(8n–5)/6 |
| **9** | 1 | 13 | 55 | 155 | 350 | 686 | 1218 | 2010 | 3135 | 4675 | n(n+1)(3n–2)/2 |
| **10** | 1 | 14 | 60 | 170 | 385 | 756 | 1344 | 2220 | 3465 | 5170 | n(n+1)(10n–7)/6 |
| **(d)** | 0 | 1 | 5 | 15 | 35 | 70 | 126 | 210 | 330 | 495 | |

## V. MULTIDIMENSIONAL SOLID NUMBERS

It seems that the Pythagorean conception of number as collections of units distributed in two- and three-dimensional spaces according to geometric form lends itself naturally to the conceptualization of multidimensional number: Distinct number entities in the same space compile to produce number formations in a higher dimensional space. This process appears endless, generating more and more complex number formations and even higher dimensional spaces where new compiling phases take place.

Let $S(v,d,n)$ denote the *n*th hypersolid number in *v* dimensions generated from an arithmetical progression with the first term *1* and common difference *d*, after *v–1* compilations. Thus, for $v = 2, 3,$ and *4,* we have

$$S(2,d,n) = p(d,n) = \sum_{r=1}^{n} x_r = \frac{1}{2} n[2+(n-1)d],$$

$$S(3,d,n) = P(d,n) = \sum_{r=1}^{n} p(d,r) = \frac{1}{2\cdot 3} n(n+1)[3+(n-1)d],$$

$$S(4,d,n) = \Pi(d,n) = \sum_{r=1}^{n} P(d,r) = \frac{1}{2\cdot 3\cdot 4} n(n+1)(n+2)[4+(n-1)d].$$

The emerging pattern suggests therefore that $S(v,d,n)$ has the form

$$S(v,d,n)=\frac{n(n+1)(n+2)\cdots(n+v-2)}{2\cdot 3\cdot 4\cdots v}[v+(n-1)d].$$

Expressed in terms of the binomial coefficients

$$\binom{M}{m}=\frac{M!}{m!(M-m)!}$$

for $m,M$ nonnegative integers $(m\le M)$, $S(v,d,n)$ reads

$$S(v,d,n)=\binom{v+n-2}{v-1}+d\cdot\binom{v+n-2}{v}. \qquad (12)$$

Next, it must be shown that the *S*-numbers defined by (12) classify as hypersolids, i.e., they result from the compiling/summation of *S*-numbers in lower dimensional spaces. Thus, it must be shown that

$$S(v,d,n)=\sum_{r=1}^{n}S(v-1,d,r), \qquad (13)$$

for any $v$ and $d$.

Indeed, using (12) for $v\to v-1$ and the binomial identity

$$\binom{m}{m}+\binom{m+1}{m}+\cdots+\binom{M}{m}=\binom{M+1}{m+1}, \qquad (14)$$

the right-hand side of (13) yields

$$\sum_{r=1}^{n}S(v-1,d,r)=\sum_{r=1}^{n}\binom{r+v-3}{v-2}+d\cdot\sum_{r=1}^{n}\binom{r+v-3}{v-1}$$

$$=\binom{v+n-2}{v-1}+d\cdot\binom{v+n-2}{v}$$

$$=S(v,d,n),$$

as required.

It should be noted that (12) is also valid for $v=0$ and $v=1$. Thus, for $v=0$ we have

$$S(0,d,n)=\begin{cases}0 & \text{if } n=0,1,\\ d & \text{if } n\ge 2,\end{cases} \qquad (15)$$

i.e., the constant differences of all arithmetical progressions (with at least two terms) may be said to constitute the zero-dimensional numbers. Also, for $v=1$ we have

$$S(1,d,n)=\begin{cases}0 & \text{if } n=0,\\ 1+d(n-1) & \text{if } n\ge 1,\end{cases} \qquad (16)$$

i.e., the individual terms of all arithmetical progressions may be said to constitute the one-dimensional numbers.

For $d=0$, (12) gives

$$S(v,0,n) = \begin{cases} 0 & \text{if } v,\ n = 0, \\ \binom{v+n-2}{v-1} & \text{if } v,\ n \geq 1. \end{cases} \tag{17}$$

Finally, for $n=0$ and any $v,\ d$, (12) reduces to

$$S(v,d,0) = 0, \tag{18}$$

while for $n = 1,$ it becomes

$$S(v,d,1) = \begin{cases} 0 & \text{if } v = 0, \\ 1 & \text{if } v \geq 1, \end{cases} \tag{19}$$

i.e., the monad is the first term of all $S$-number sequences. Subsequently, all hypersolid numbers are, strictly speaking, pyramidal in form.

To summarize, multidimensional solid numbers, or hypersolids, in short, are defined by relations (12) and (13) for any $v,d,n \geq 0$ .

It readily follows from (13) that

$$S(v,d,n) = S(v,d,n-1) + S(v-1,d,n), \tag{20}$$

i.e., the $n$-gnomons of the hypersolids $S(v,d,n-1)$ are the hypersolids $S(v-1,d,n)$ defined in the immediately lower dimensional space. Also, the $d$-gnomons for any given $v$-dimensional space are the first order $(d = 1)$ hypersolids:

$$S(v,d,n) = S(v,d-1,n) + S(v,1,n-1). \tag{21}$$

## VI. DIMENSIONAL GNOMONS AND ARITHMETIC TRIANGLES

Hypersolids in consecutive dimensional spaces with the same $n$ and $d$ differ by what may be defined as dimensional or $v$-gnomons. Thus, given the sequences $S(v-1,d,1)$, $S(v-1,d,2)$,..., and $S(v,d,1)$, $S(v,d,2)$,..., $v$-gnomons are the differences of their corresponding terms:

$$S(v,d,n) - S(v-1,d,n).$$

By use of (12) one finds that $v$-gnomons are the hypersolids $S(v,d,n-1)$. Hence, relation (20) proves to be gnomonic in a two-fold sense, i.e., it defines symmetrically both the $n$-gnomons and the $v$-gnomons.
The hypersolid numbers generated from the same arithmetical progression by successive compilations may best be displayed in a triangular form, as shown in Table 4.

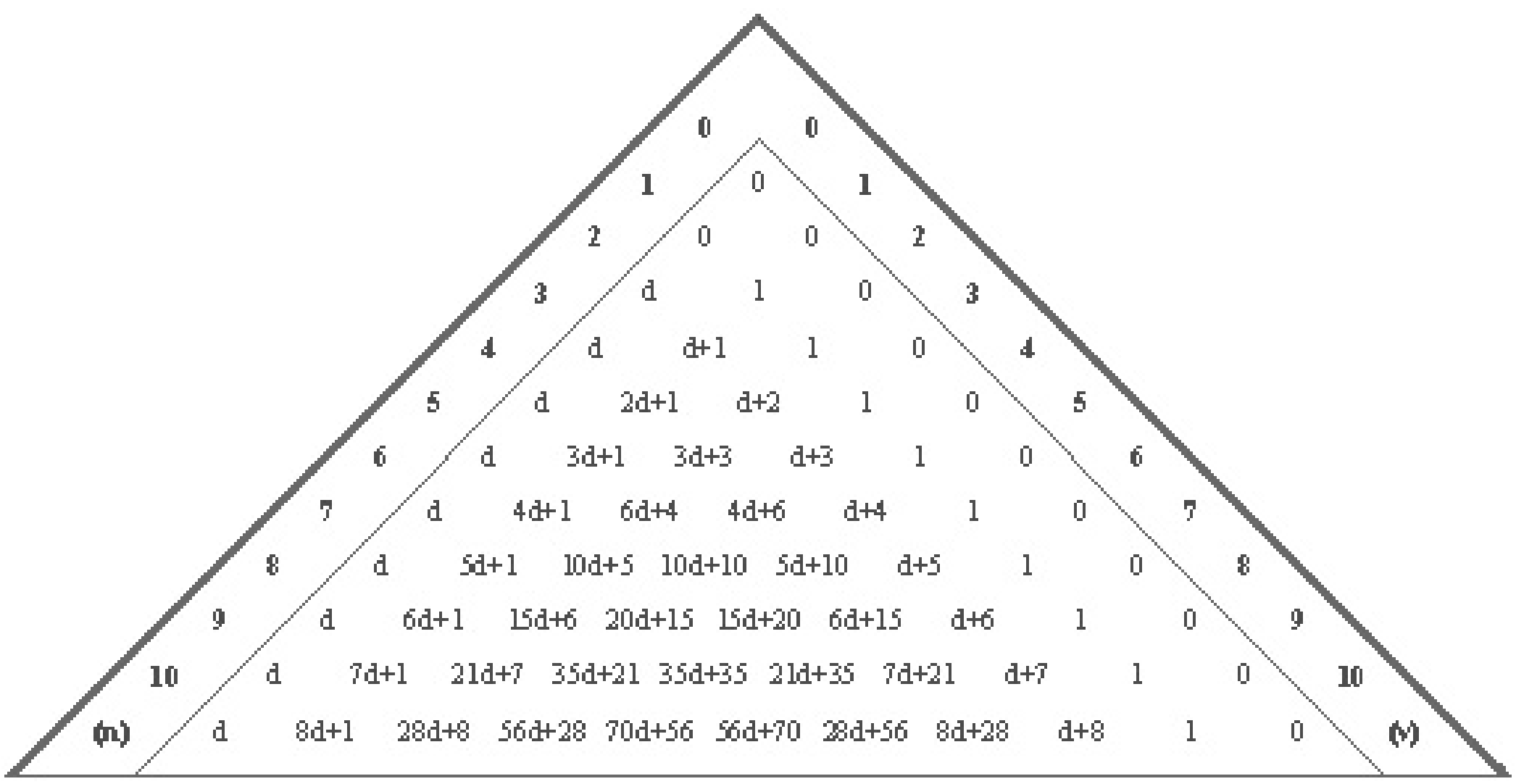

**Table 4 Hypersolid Numbers with the Same d**

From the resulting pattern a number of regularities becomes manifest of which the principal ones seem to be the following:

**Corollary 1.** *Each entry in this arithmetic triangle with* $n+v \geq 3$ *is equal to the sum of its two adjacent entries in the horizontal row above.* Evidently, this sum rule is relation (20) at work.

**Corollary 2.** *The sum of the first n entries along any v-row is equal to the nth entry in the next* $(v+1)$*-row.* This manifests directly the compiling process of generating the hypersolids as defined by (13).

**Corollary 3.** *The sum of the first v entries along any n-row is equal to the vth entry in the next* $(n+1)$*-row.* Thus, we have

$$S(v,d,n+1) = \sum_{r=0}^{v} S(r,d,n), \tag{22}$$

for fixed *d* and *v*. Indeed, the right-hand side of (22) can be written as

$$\sum_{r=0}^{v} S(r,d,n) = \sum_{r=1}^{v} \binom{n+r-2}{r-1} + d \cdot \sum_{r=0}^{v} \binom{n+r-2}{r}$$

which, by virtue of the binomial identity

$$\binom{M}{0} + \binom{M+1}{1} + \cdots + \binom{M+m}{m} = \binom{M+m+1}{m}, \tag{23}$$

yields

$$\binom{v+n-1}{v-1} + d \cdot \binom{v+n-1}{v} = S(v,d,n+1).$$

**Corollary 4.** *Along each horizontal row, the sum of all entries equals*

$$\sum_{v+n=c} S(v,d,n) = \begin{cases} 0 & \text{if } c < 2, \\ (d+1) \cdot 2^{c-2} & \text{if } c \geq 2, \end{cases} \tag{24}$$

*with* $c = v + n$ *a constant*. This follows directly after expanding the left-hand side and applying the fundamental binomial identity

$$\binom{M}{0} + \binom{M}{1} + \cdots + \binom{M}{M} = 2^M. \tag{25}$$

It must be noted that for $v+n \geq 3$ the sum of the entries in each horizontal row is twice that of the preceding horizontal row and, if diminished by $d+1$, is also equal to the sum of the entries in all the preceding horizontal rows.

**Corollary 5.** *Summing the entries lying on directional lines of slope 1/2 , one obtains the Fibonacci sequences*

$$d, \quad d+1, \quad 2d+1, \quad 3d+2, \quad 5d+3, \quad 8d+5, \quad 13d+8, \quad 21d+13, \quad 34d+21, \ldots,$$

which satisfy the recurrence relation of second order $a_n = a_{n-1} + a_{n-2}$ given that $a_0 = d$ and $a_1 = d+1$. Other sequences satisfying recurrence relations of higher order may be obtained by summing entries along lines of smaller slope.

**Corollary 6.** *For d = 0, the arithmetic triangle reduces to that of Pascal.* Thus, Pascal´s triangle manifests itself to be a particular cross section of the universe of hypersolid numbers.

## VII. SUMS OF HYPERSOLID NUMBERS

Having established that hypersolid numbers are characterized by the three parameters *v, d* and *n,* in what follows we shall determine the sums of hypersolid numbers with fixed $s = v + d + n$ and a given *v, d* and *n*, respectively. Subsequently, we shall determine the sum of all hypersolid numbers with the same *s*. It is clear from relations (15), (16) and (17) that non-zero sums can be obtained for $s \geq 2$ .

**Theorem 1.** *The sum of (non-zero) hypersolid numbers with same s and v, equals*

$$\sum_{s,v\ given} S(v,d,n) = \begin{cases} \binom{s-1}{2} & if\ v = 0, \\ \binom{s-1}{v} + \binom{s-1}{v+2} & if\ v \geq 1, \end{cases} \tag{26}$$

*and their multitude is* $s-2$ *for* $v = 0$, *and* $s-v$ *for* $v \geq 1$.

From the well-known fact that the number of non-negative solutions to the equation $x_1 + x_2 + ... + x_m = M$ is

$$\binom{M+m-1}{M},$$

it trivially follows that the multitude of all hypersolid numbers with the same *s* and *v* is $s-v+1$. For $v = 0$ , there are three zero *S*-numbers to be discarded, namely, the numbers $S(0,s,0)$, $S(0,s-1,1)$ and $S(0,0,s)$. Hence, the non-zero *S*-numbers in this case are $s-2$. For $v \geq 1$, there is only one zero *S*-number, namely, the number $S(v,s-v,0)$ , and the multitude of non-zero *S*-numbers in this case is $s-v$.

For $v = 0$, the left-hand side of (26) expanded reads

$$\sum_{d=0}^{s} d \cdot \binom{s-d-2}{0} = \sum_{d=1}^{s-2} d = \binom{s-1}{2},$$

as required. For $v \geq 1$, we have

$$\sum_{d=0}^{s-v-1} \binom{s-d-2}{v-1} + \sum_{d=0}^{s-v-2} d \cdot \binom{s-d-2}{v}.$$

Using (23) and the identity

$$1 \cdot \binom{M-1}{m} + 2 \cdot \binom{M-2}{m} + 3 \cdot \binom{M-3}{m} + \cdots + (M-m) \cdot \binom{m}{m} = \binom{M+1}{m+2} \tag{27}$$

with $M > m$ , we get

$$\binom{s-1}{v} + \binom{s-1}{v+2}.$$

The proof of Theorem 1 is thus complete.

As an illustration of Theorem 1, some indicative sums of *S*-numbers from Tables 1, 2 and 3 are actually performed in Table 5 below, thus verifying the answers readily derived from (26).

**Table 5 Sums of Hypersolid Numbers with *s, v* Given**

| s | v | Actual Sums | Theorem 1 |
|---|---|---|---|
| **10** | **2** | 1 + 8 + 18 + 28 +35 +36 + 28 + 8 | $\binom{9}{2} + \binom{9}{4} = 162$ |
| **11** | **2** | 1 + 9 + 21 + 34 + 45 + 51 + 49 + 36 + 9 | $\binom{10}{2} + \binom{10}{4} = 255$ |
| **12** | **3** | 1 + 10 + 30 + 60 + 95 + 126 + 140 + 120 + 45 | $\binom{11}{3} + \binom{11}{5} = 627$ |
| **13** | **3** | 1 + 11 + 34 + 70 + 115 + 161 + 196 + 204 + 165 + 55 | $\binom{12}{3} + \binom{12}{5} = 1012$ |
| **14** | **4** | 1 + 12 + 45 + 110 + 210 + 336 + 462 + 540 + 495 + 220 | $\binom{13}{4} + \binom{13}{6} = 2431$ |
| **15** | **4** | 1 + 13 + 50 + 125 + 245 + 406 + 588 + 750 + 825 + 715 + 286 | $\binom{14}{4} + \binom{14}{6} = 4004$ |

The sums in Table 5 are complemented by the corresponding $d = 0$ numbers not listed in these Tables.

**Theorem 2.** *The sum of (non-zero) hypersolid numbers with same s and d equals*

$$\sum_{s,d\ given} S(v,d,n) = (d+1)\cdot 2^{s-d-2} \tag{28}$$

*for $s \geq 2$, and their multitude is $s-1$ for $d = 0$, and $s - d$ for $d \geq 1$.*

Evidently, this Theorem restates (24) for $c = v+n = s-d$. The same argument used in Theorem 1 implies that all possible hypersolid numbers with the same $s$ and $d$ are $s-d+1$. For $d = 0$, there are two zero $S$-numbers to be discarded, namely, the numbers $S(s,0,0)$, $S(0,0,s)$ and, consequently, the multitude of non-zero $S$-numbers is $s-1$. For $d \geq 1$, only number $S(s-d,d,0)$ is zero; hence, the multitude of non-zero $S$-numbers in this case is $s-d$.

A most representative illustration of (28) is provided in Table 4 (p. 8).

**Theorem 3.** *The sum of (non-zero) hypersolid numbers with same s and n equals*

$$\sum_{s,n\ given} S(v,d,n) = \begin{cases} s-1 & if\ n=1, \\ 2\cdot\binom{s-1}{n} & if\ 2 \leq n < s, \end{cases} \tag{29}$$

*and their multitude is $s-n$ for $n=1$ and $s-n+1$ for $n \geq 2$.*

By the same argument employed in Theorems 1 and 2, all hypersolid numbers with the same $s$ and $n$ are $s-n+1$. For $n \geq 1$, all these numbers are non-zero except the number $S(0,s-1,1)$. Consequently, the multitude of non-zero $S$-numbers is $s-n$ for $n=1$, and $s-n+1$ for $n \geq 2$.

Expanding the left-hand side of (29), we have

$$\sum_{v=0}^{s-n}\binom{v+n-2}{v-1} + \sum_{v=0}^{s-n}(s-n-v)\cdot\binom{v+n-2}{v}.$$

For $n = 0$, the binomial coefficients are zero and both sums are zero. For $n = 1$, the second sum is zero while the first sum equals $s-1$. Using (14) for n$\geq 2$, we get

$$= \binom{s-1}{n} + (s-n)\cdot\binom{s-1}{n-1} - \sum_{v=0}^{s-n} v\cdot\binom{v+n-2}{v}$$

$$= \binom{s-1}{n} + (s-n)\cdot\binom{s-1}{n-1} - \frac{1}{(n-2)!}\sum_{v=1}^{s-n} P_{v-1}^{v+n-2},$$

with $P_m^M$ the usual permutation coefficients. Then, applying the identity

$$P_m^{M+1} = m! + m\cdot(P_{m-1}^{M} + P_{m-1}^{M-1} + \cdots + P_{m-1}^{m}) \tag{30}$$

with $m \leq M$, one finds

$$\binom{s-1}{n} + \frac{(s-1)!}{(n-2)!(s-n-1)!}\left(\frac{1}{n-1} - \frac{1}{n}\right) = 2\cdot\binom{s-1}{n},$$

as desired.

As an illustration of Theorem 3, Table 6 displays a particular cross section of hypersolid numbers with *n* fixed as well as the relevant sums predicted by (29).

**Table 6 Sums of Hypersolids with *s*, *n* Given (*n=4; d=0, ···, 10*)**

| $2\cdot\binom{s-1}{4}$ | | 2·0 ↑ | 2·1 ↑ | 2·5 ↑ | 2·15 ↑ | 2·35 ↑ | 2·70 ↑ | 2·126 ↑ | 2·210 ↑ | 2·330 ↑ | 2·495 ↑ | 2·715 ↑ |
|---|---|---|---|---|---|---|---|---|---|---|---|---|
| **s** | **d** \ ν | **0** | **1** | **2** | **3** | **4** | **5** | **6** | **7** | **8** | **9** | **10** |
| **4** | **0** | 0 | 1 | 4 | 10 | 20 | 35 | 56 | 84 | 120 | 165 | 220 |
| **5** | **1** | 1 | 4 | 10 | 20 | 35 | 56 | 84 | 120 | 165 | 220 | 286 |
| **6** | **2** | 2 | 7 | 16 | 30 | 50 | 77 | 112 | 156 | 210 | 275 | 352 |
| **7** | **3** | 3 | 10 | 22 | 40 | 65 | 98 | 140 | 192 | 255 | 330 | 418 |
| **8** | **4** | 4 | 13 | 28 | 50 | 80 | 119 | 168 | 228 | 300 | 385 | 484 |
| **9** | **5** | 5 | 16 | 34 | 60 | 95 | 140 | 196 | 264 | 345 | 440 | 550 |
| **10** | **6** | 6 | 19 | 40 | 70 | 110 | 161 | 224 | 300 | 390 | 495 | 616 |
| **11** | **7** | 7 | 22 | 46 | 80 | 125 | 182 | 252 | 336 | 435 | 550 | 682 |
| **12** | **8** | 8 | 25 | 52 | 90 | 140 | 203 | 280 | 372 | 480 | 605 | 748 |
| **13** | **9** | 9 | 28 | 58 | 100 | 155 | 224 | 308 | 308 | 525 | 660 | 814 |
| **14** | **10** | 10 | 31 | 64 | 110 | 170 | 245 | 336 | 444 | 570 | 715 | 880 |

By applying any of the above Theorems, we can now find the sum of hypersolid numbers with the same *s*. Thus, we have the following Theorem:

**Theorem 4.** *The sum of (non- zero) hypersolid numbers with the same* $s = v + d + n$*, equals*

$$\sum_{s\ given} S(v,d,n) = 2^s - (s+1), \qquad (31)$$

*and their multitude is* $\binom{s+1}{2} - 2,\ \ for\ \ s \geq 2.$

All *S*-numbers with $v + d + n = s$ fixed, are $\binom{s+2}{2}$. Among these, there are $s+1$ zero numbers for $n = 0$, one for $v = d = 0$, and one for $v = 0$ and $n = 1$. Consequently, the multitude of non-zero *S*-numbers is

$$\binom{s+2}{2} - (s+1) - 2 = \binom{s+1}{2} - 2.$$

To show (31), we shall use Theorems 1, 2, and 3 in turn:

- By Theorem 1, the sum in (31) may be written as

$$\binom{s-1}{2} + \sum_{v=1}^{s} \binom{s-1}{v} + \sum_{v=1}^{s} \binom{s-1}{v+2}$$

which, using (25), becomes

$$= \binom{s-1}{2} + (2^{s-1} - 1) + [\ 2^{s-1} - \binom{s-1}{0} - \binom{s-1}{1} - \binom{s-1}{2}\ ]$$

$$= 2 \cdot 2^{s-1} - \binom{s-1}{0} - \binom{s-1}{1} = 2^s - (s+1).$$

- By Theorem 2, the sum in (31) for $r = s - d - 2$ may be written as

$$\sum_{d=0}^{s-2} (d+1) \cdot 2^{s-d-2} = (s-1) \cdot \sum_{r=0}^{s-2} 2^r - \sum_{r=0}^{s-2} r \cdot 2^r ,$$

and using the identities

$$\sum_{r=0}^{R} 2^r = 2^{R+1} - 1 , \quad \sum_{r=0}^{R} r \cdot 2^r = 2 + (R-1) \cdot 2^{R+1} , \qquad (32)$$

it becomes

$$(s-1) \cdot (2^{s-1} - 1) - (s-3) \cdot 2^{s-1} - 2 = 2^s - (s+1),$$

as required.

- By Theorem 3, the sum in (31) may be written as

$$(s-1) + 2 \cdot \sum_{n=2}^{s-1} \binom{s-1}{n} ,$$

and using (25), one obtains

$$(s-1) + 2 \cdot (2^{s-1} - s) = 2^s - (s+1),$$

as required.

As an illustration of Theorem 4, we list below the non-zero hypersolid numbers with $s=6$:

| | | | | |
|---|---|---|---|---|
| $S(0,1,5) = 1,$ | $S(0,2,4) = 2,$ | $S(0,3,3) = 3,$ | $S(0,4,2) = 4,$ | |
| $S(1,0,5) = 1,$ | $S(1,1,4) = 4,$ | $S(1,2,3) = 5,$ | $S(1,3,2) = 4,$ | $S(1,4,1) = 1,$ |
| $S(2,0,4) = 4,$ | $S(2,1,3) = 6,$ | $S(2,2,2) = 4,$ | $S(2,3,1) = 1,$ | |
| $S(3,0,3) = 6,$ | $S(3,1,2) = 4,$ | $S(3,2,1) = 1,$ | | |
| $S(4,0,2) = 4,$ | $S(4,1,1) = 1,$ | | | |
| $S(5,0,1) = 1.$ | | | | |

Therefore, we have $\binom{7}{2} - 2 = 19$ $S$-numbers, and their sum equals $2^6 - 7 = 57$.

# Appendix A

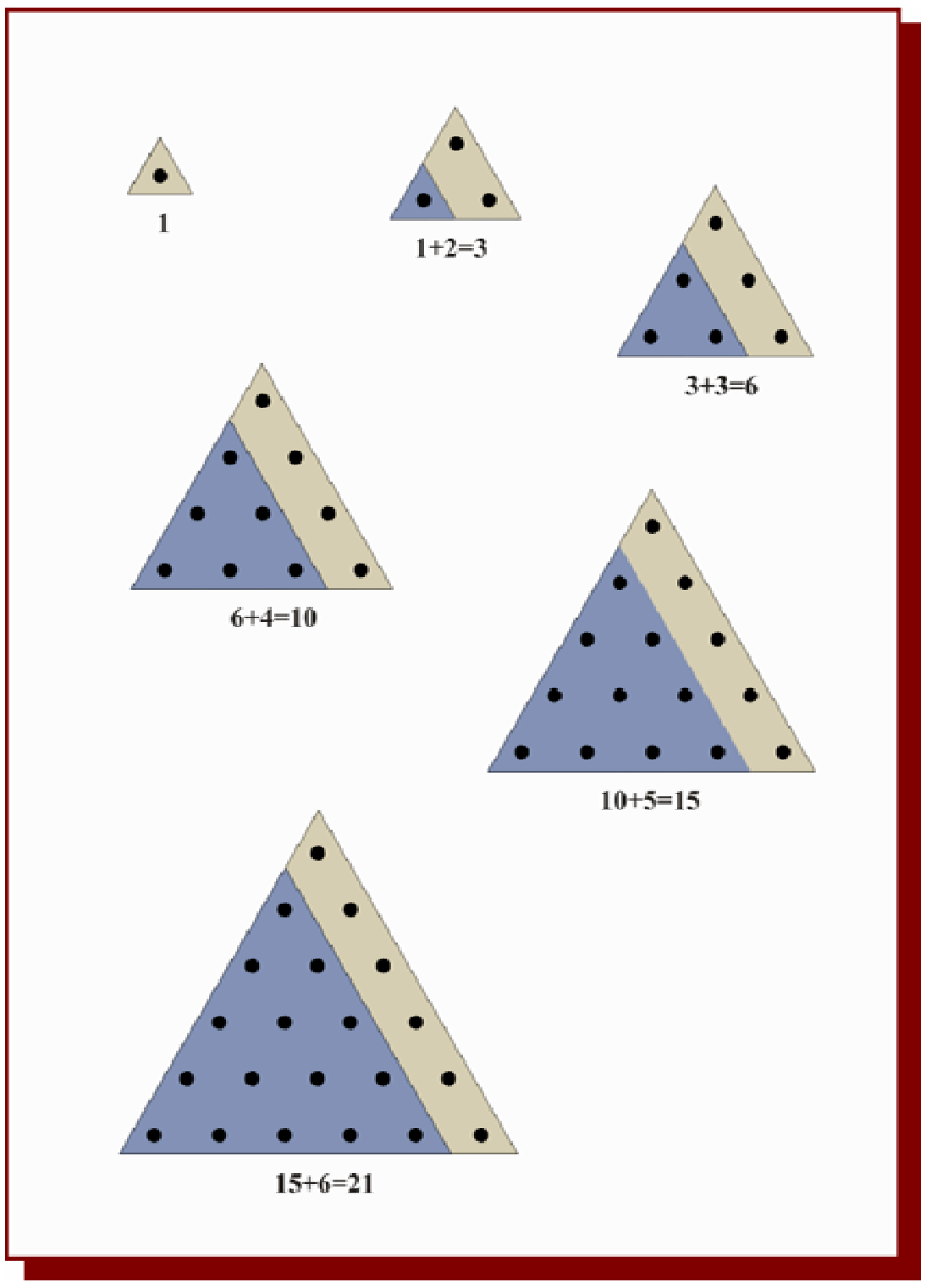
1
1+2=3
3+3=6
6+4=10
10+5=15
15+6=21

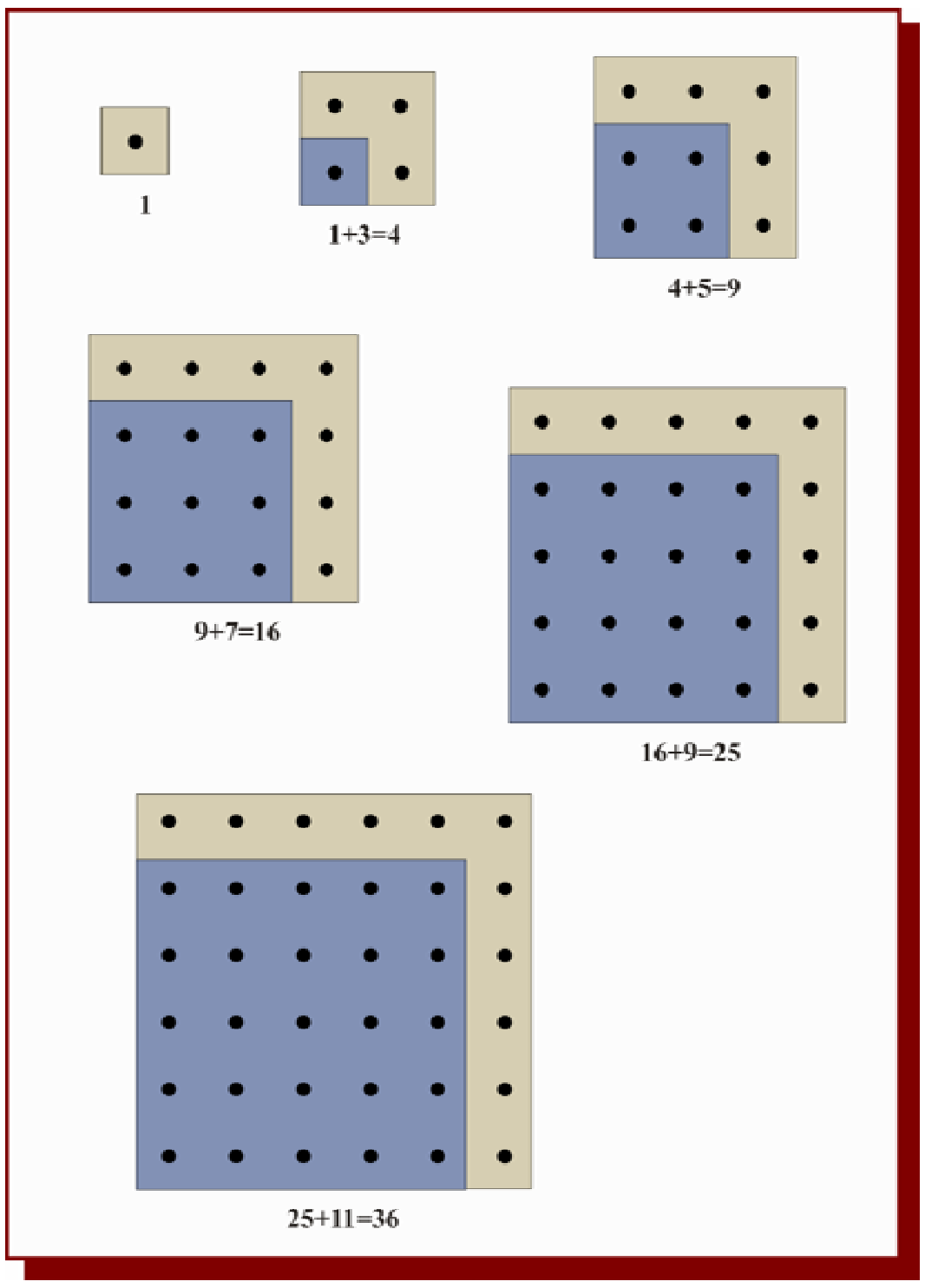
1
1+3=4
4+5=9
9+7=16
16+9=25
25+11=36

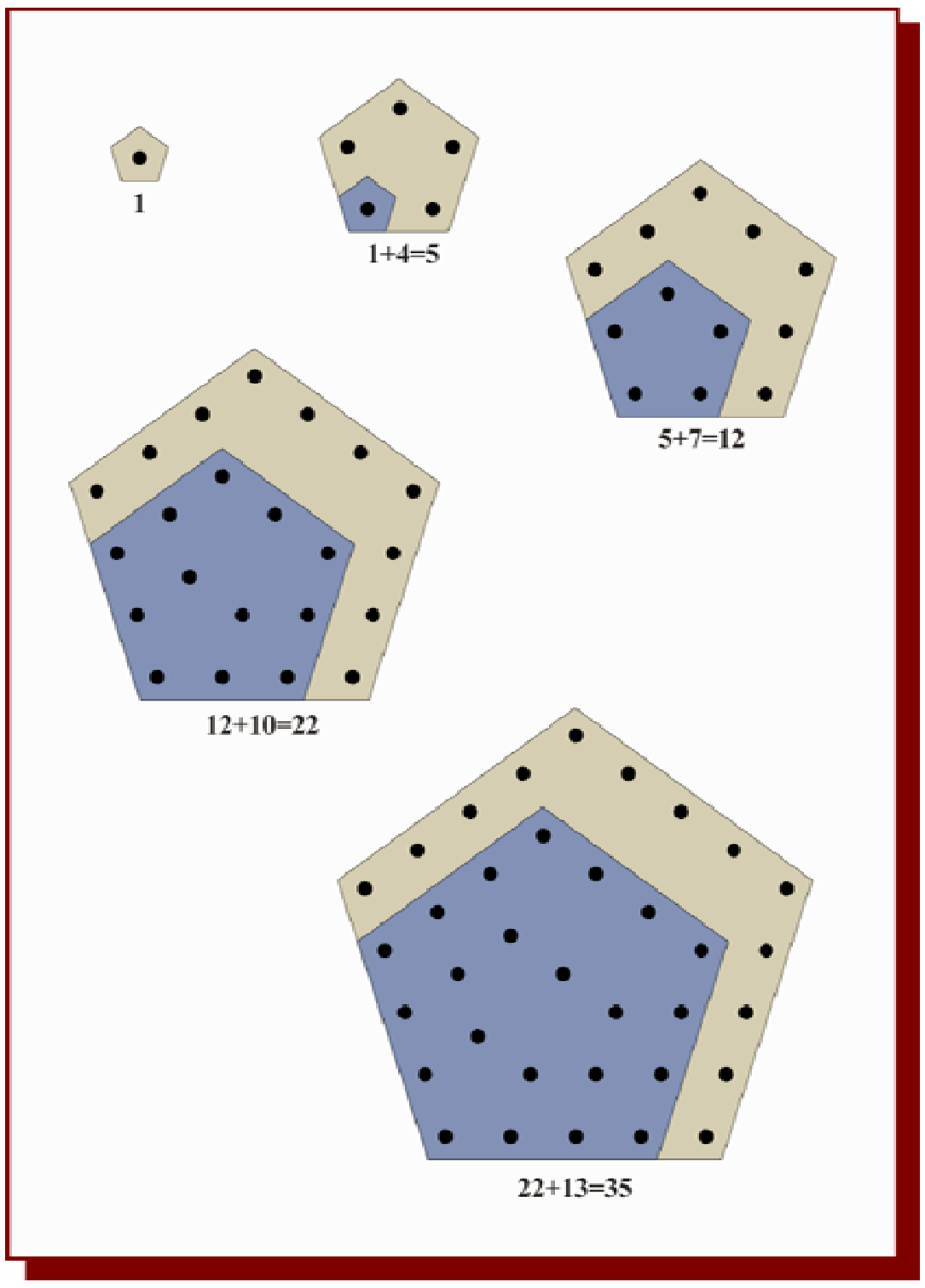
1
1+4=5
5+7=12
12+10=22
22+13=35

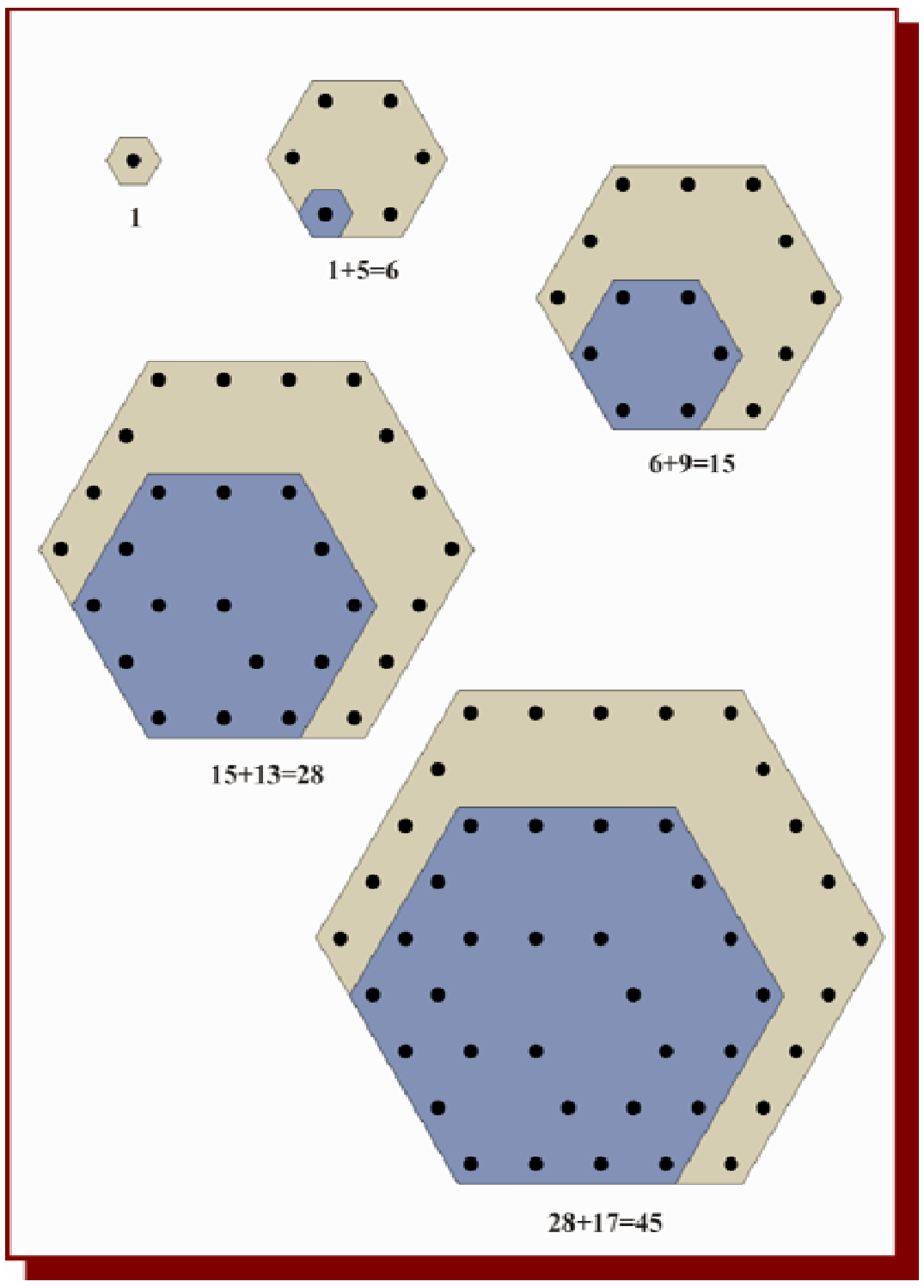
1
1+5=6
6+9=15
15+13=28
28+17=45

Appendix B

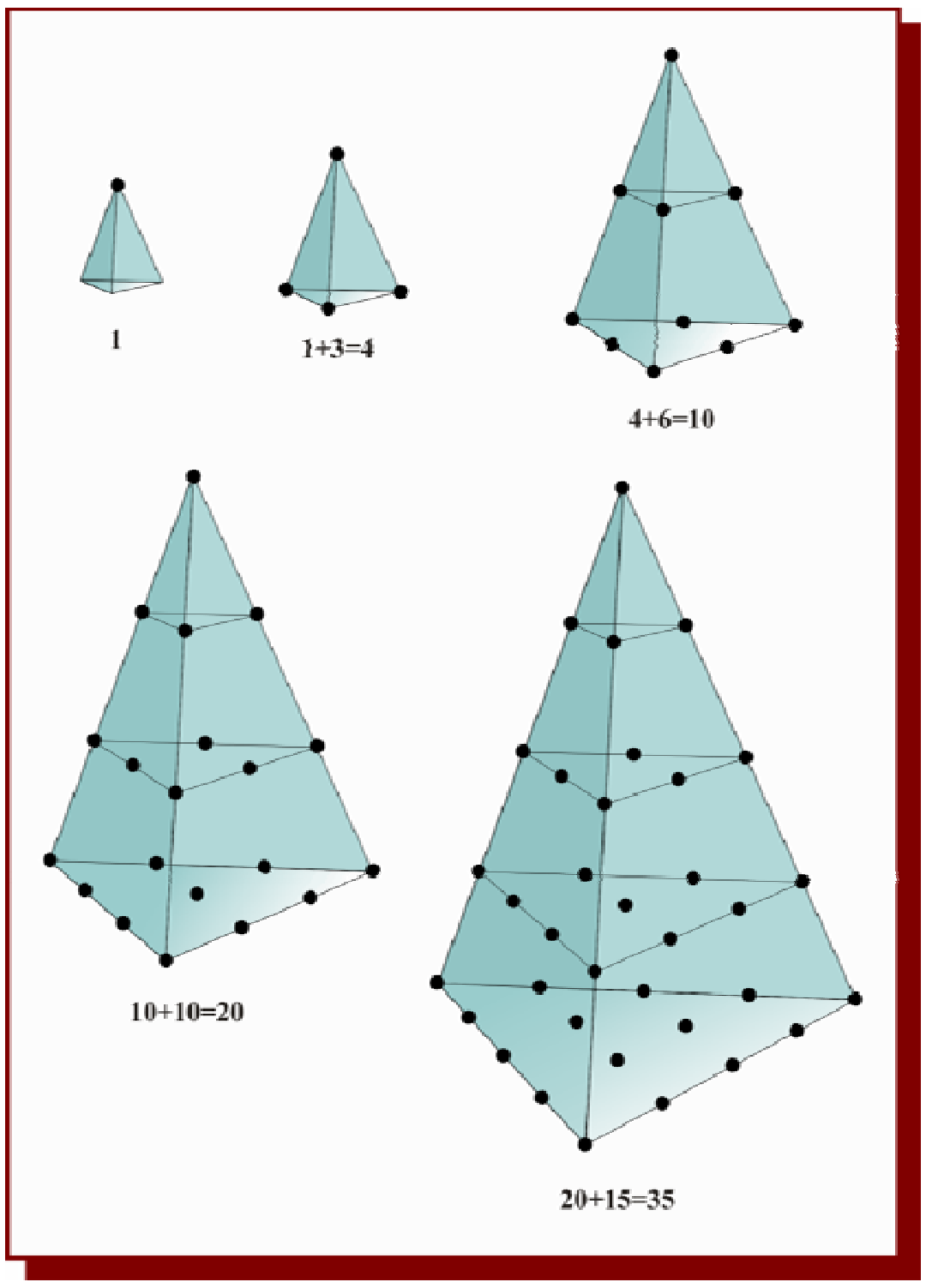

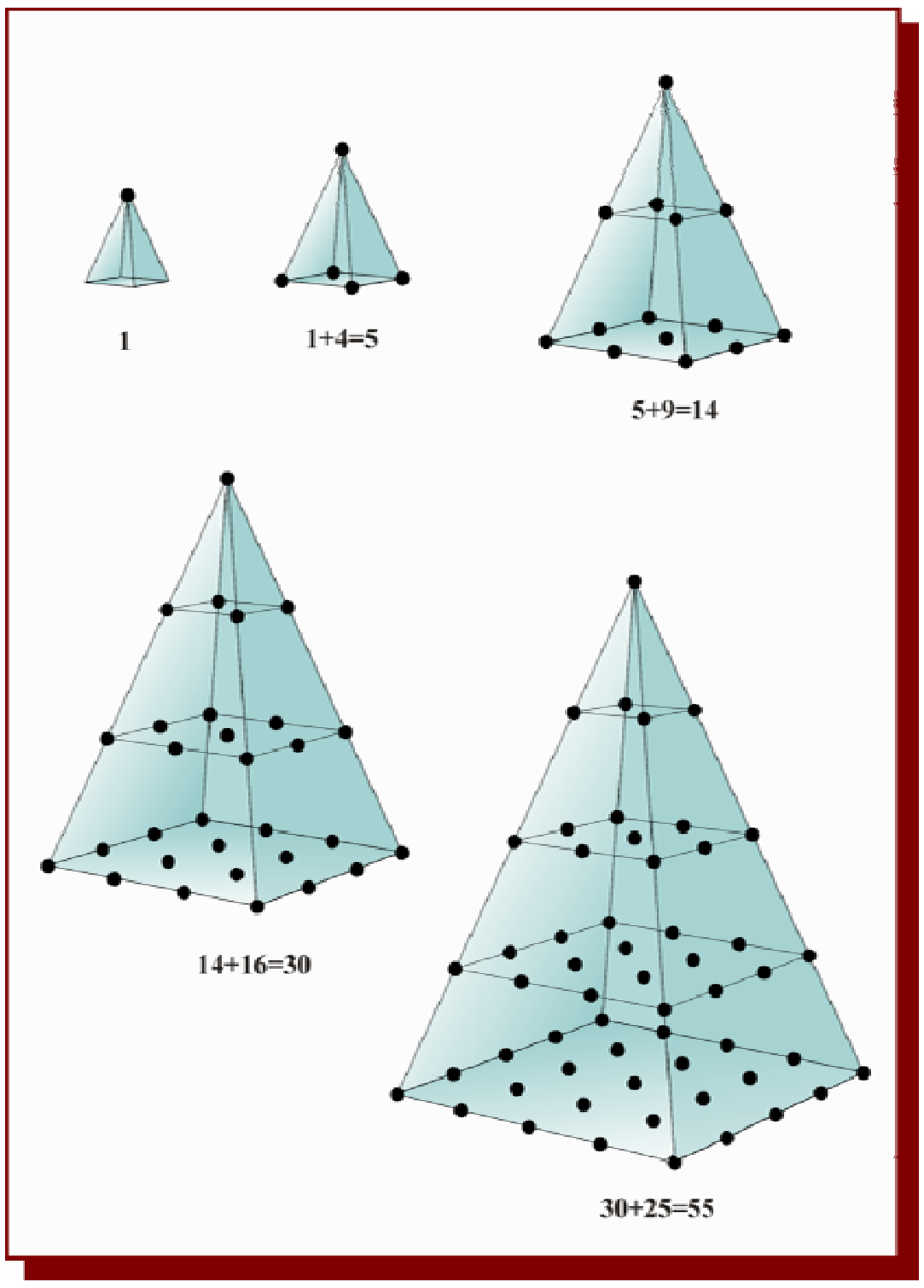